\documentclass[onefignum,onetabnum]{siamart190516}

\usepackage[margin=1in]{geometry}
\usepackage{amsmath}
\usepackage{amssymb}
\usepackage{graphicx}
\usepackage{epstopdf}
\usepackage{onimage}
\usepackage[mathcal]{euscript}

\graphicspath{{Figures/}}

\DeclareMathOperator{\sech}{sech}
\newcommand{\J}{J}
\newcommand{\scriptJ}{\mathcal{J}}
\newcommand{\AR}{\text{AR}}
\newcommand{\MR}{\text{MR}}

\title{Optimal control with learning on the fly: a toy problem}
\date{26 February 2020}
\author{C. L. Fefferman \and B. Guill\'en Pegueroles \and C. W. Rowley \and M. Weber}
\begin{document}
\maketitle

\begin{abstract}
  We exhibit optimal control strategies for a simple toy problem in which the
  underlying dynamics depend on a parameter that is initially unknown and must
  be learned.
  We consider a cost function posed over a finite time interval, in contrast to
  much previous work that considers asymptotics as the time horizon tends to infinity.
  We study several different versions of the problem, including
  Bayesian control, in which we assume a prior distribution on the unknown
  parameter; and ``agnostic'' control, in which we assume nothing about the
  unknown parameter.  For the agnostic problems, we compare our performance with that of an opponent who knows the value of the
  parameter. This comparison gives rise to several notions of ``regret,'' and we
  obtain strategies that minimize the ``worst-case regret'' arising from the most
  unfavorable choice of the unknown parameter.  In every case, the optimal
  strategy turns out to be a Bayesian strategy or a limit of Bayesian
  strategies.
\end{abstract}



\section{Motivation and introduction}
\label{sec:intro}

We investigate control problems in which we must make decisions with little time
and little data available.  Our motivating example is the success of pilots
learning in real time to fly and safely land an airplane after it has been
severely damaged, for instance as documented in \cite{Brazy:2009}.

Control with learning has been considered in many different application areas;
see, for example, books such as \cite{Cesa-Bianchi:2006,Hazan:2016,Bertsekas:2005,Powell:2011} as well as
recent papers such as \cite{Agarwal:2019,Agarwal:2019a}.  Much is known also
about the closely related ``multi-armed bandit'' problem; see, for instance, the
classic papers \cite{Lai:1985,Auer:2002} and the more recent
survey~\cite{Bubeck:2012}.

A standard approach to control with learning is to start by estimating the
parameters in the model, and then to design a controller that is optimal for
that model.  A related approach is to divide the available time into epochs, and
within each epoch, first refine the model, and then update the controller
accordingly.  In~\cite{Cohen:2019}, it is shown that this latter approach gives
results that are optimal (in some sense) in the asymptotic limit of large time.
For the problem we consider, such an approach cannot lead to an optimal control
strategy.  We are interested in results that are optimal for a
fixed time interval, rather than for large time.  Our optimal strategies do not
divide into distinct phases of exploration and exploitation, but instead take
past history into account at every moment.


This article analyzes a toy problem in which we apply a time-dependent control
to keep the position of a moving particle close to zero.  The position
$q(t)\in\mathbb{R}$ is governed by Brownian motion with a drift.  The drift rate
is given by $a + u(t)$, where $a$ is an unknown constant, and $u(t)$ is our
control at time~$t$.  For this simple one-dimensional toy model, we consider
several different notions of optimality, and we exhibit control policies that
are optimal with respect to each of these.  At the end of the paper, we mention
two additional toy problems to which we hope to return in later work.

\section{The toy problem}
\label{sec:toy-problem}

We begin with a time interval $[0,T]$ subdivided into discrete timesteps of size~$\Delta t$.  From $t$ to $t+\Delta t$, the change in the position $q$ is given
by
\begin{equation}
  \label{eq:1}
  \Delta q = (a + u(t)) \Delta t + \Delta W,
\end{equation}
where $a$ is an unknown constant, $u(t)$ is our control, and $\Delta W$ is a
normally distributed random variable with zero mean and standard deviation
$\sigma_0 (\Delta t)^{1/2}$ for a known coefficient~$\sigma_0$.  A simple
scaling allows us to take $\sigma_0=1$, which we do from now on.  In the limit as
$\Delta t$ tends to zero, we obtain a control system in continuous time, for
which $W(t)$ is Brownian motion and $dW(t)$ is white noise.  Our goal is to find
optimal control strategies for this continuous-time system.



From time 0 to some given time $T_0$, we are allowed only to observe the
particle: i.e., we must take $u=0$.  From time $T_0$ to time $T$, we may apply
any control strategy we please, provided that $u(t)$ is determined by the
history $q(\tau)$ for $\tau\le t$.  We want to pick a control strategy to
minimize the expected value of a cost function
\begin{equation}
  \label{eq:32}
  \J = \int_{T_0}^T (q^2 + \lambda u^2)\,dt,
\end{equation}
where $\lambda>0$ is a known coefficient.  If the parameter~$a$ is known, then
the only randomness in the system arises from the Brownian motion, so the notion
of expected value is well defined.  If $a$ is unknown, the meaning of the
expected value is not immediately clear.  We will discuss it carefully in
Sections \ref{sec:bayesian-control} and~\ref{sec:agnostic-control}.

The most interesting case arises when $T_0=0$. Also, a simple scaling allows us
to take $\lambda=1$ without disturbing the normalization $\sigma_0=1$.  Unless we
say otherwise, we will assume that $\lambda=1$ and~$T_0=0$.  Furthermore, we
suppose that our particle starts at position $q(0)=0$.

\section{Notions of optimality}
\label{sec:optimality}

In this section, we provide careful definitions of optimal strategies, first
assuming that the parameter~$a$ is known (classical control), then assuming a
prior belief regarding~$a$ (Bayesian control), and finally assuming no prior
knowledge of~$a$ (``agnostic'' control).

\subsection{Known parameter value~$a$}
\label{sec:known-param}

If the value of~$a$ is known, then we simply ask for a control strategy that
minimizes the expected value of~(\ref{eq:32}).  As observed above, the expected
value is well defined because $a$ is known.  To calculate such an optimal
control strategy is a classical problem of control theory, and we review its
solution in Section~\ref{sec:classical-optimal} below.

\subsection{Bayesian control}
\label{sec:bayesian-control}

Next, suppose that we are given a probability measure $\mu(a)$
reflecting our prior belief about the unknown~$a$.  That is, the probability
that $a$ lies between $\alpha$ and $\beta$ is given by
\begin{equation*}
  \int_\alpha^\beta d\mu(a).
\end{equation*}
Then as in the classical case, we can make sense of the expected value of~$\J$.
We compute the expected value of~$\J$ assuming a given value of $a$; call that
quantity $\J(a)$.  We then average over all~$a$ according to the probability measure~$\mu$:
\begin{equation*}
  E(\J) = \int \J(a)\,d\mu(a).
\end{equation*}
The goal of Bayesian control is to find a control strategy $u$ that minimizes
$E(\J)$.  This strategy will of course depend on our prior belief~$\mu$, and the
optimal strategy is discussed in Section~\ref{sec:bayesian-strategy}.

\subsection{Agnostic control}
\label{sec:agnostic-control}

Finally, suppose we know nothing about the parameter~$a$.  We hope to pick our
strategy to minimize one of several notions of {\em regret}, which we spell out
below.  In all variants, we play against an opponent who knows the value of~$a$
and plays perfectly, while we know nothing about~$a$.  Suppose we pick a control
strategy~$\mathcal{Q}$ without knowing~$a$.  Given the true
value of~$a$, we can compute the expected cost $\J_\mathcal{Q}(a)$ of our
strategy as in the classical case.  We want to compare $\J_\mathcal{Q}(a)$ with
the expected cost $\J_\text{opponent}(a)$ for our opponent.

We can now give precise descriptions of three problems of agnostic control:

\paragraph{Additive regret}
The {\em additive regret}, often called simply {\em regret}, is the
  difference
  \begin{equation*}
    \AR_\mathcal{Q}(a) = \J_\mathcal{Q}(a) - \J_\text{opponent}(a) \ge 0.
  \end{equation*}
  We can look for a control strategy~$\mathcal{Q}$ that minimizes the worst-case
  additive regret
  \begin{equation*}
    \AR^*_\mathcal{Q} = \sup_{a} \AR_\mathcal{Q}(a).
  \end{equation*}
\paragraph{Multiplicative regret}
  The {\em multiplicative regret}, often called the {\em competitive ratio},
  is the ratio
  \begin{equation*}
    \MR_\mathcal{Q}(a) = \frac{\J_\mathcal{Q}(a)}{\J_\text{opponent}(a)} \ge 1.
  \end{equation*}
  We can look for a control strategy~$\mathcal{Q}$ that minimizes the worst-case
  multiplicative regret
  \begin{equation*}
    \MR^*_\mathcal{Q} = \sup_a \MR_\mathcal{Q}(a).
  \end{equation*}
\paragraph{Fuel tax regret}  As before, we compute our expected cost
  using formula~(\ref{eq:32}) (which depends on our strategy~$\mathcal{Q}$),
  with $\lambda=1$.  However, we now assume that our opponent incurs a cost
  \begin{equation*}
    \J = \int_{T_0}^T (q^2 + \lambda u^2)\,dt,
  \end{equation*}
  where $\lambda>1$.  Thus, our opponent pays a {\em fuel tax}.  We
  can look for a control strategy~$\mathcal{Q}$ such that, for any assumed value
  of~$a$, our expected cost is at most that of our opponent.  We want to find
  such a strategy with $\lambda$ as small as possible.  We define the {\em
    fuel tax regret} to be the above minimal~$\lambda$.

  \bigskip
\noindent Solutions of the above agnostic control problems are shown in Section~\ref{sec:results}.

\subsection{A constant-regret Bayesian strategy is optimal}
\label{sec:constant-regret-opt}

Our solutions to the agnostic control problems will be Bayesian for a particular
choice of prior belief~$\mu$.  We explain the idea for multiplicative regret;
analogous ideas apply to the other agnostic control problems mentioned above.

Suppose that a particular prior belief $\mu$ gives rise to an optimal Bayesian
strategy $\mathcal{B}$ whose multiplicative regret $\MR_\mathcal{B}(a)$ is {\em
  constant\/} (independent of~$a$).  Then the strategy $\mathcal{B}$ minimizes
worst-case regret $\MR^*$.  To see this, we argue as follows.

Suppose that instead of~$\mathcal{B}$, we use another strategy~$\mathcal{C}$.
The strategy~$\mathcal{C}$ cannot perform better than~$\mathcal{B}$ for all
values of~$a$; otherwise $\mathcal{B}$ would not be optimal for the prior
belief~$\mu$.  Thus, there is a value of $a$ for which $\MR_\mathcal{C}(a)\ge
\MR_\mathcal{B}(a)$.  But since the strategy $\mathcal{B}$ has constant regret,
the right-hand side is independent of~$a$.  In other words, for some $a$, we
have
\begin{equation*}
  \MR_\mathcal{C}(a)\ge \MR^*_\mathcal{B}.
\end{equation*}
Consequently, the worst-case regret of $\mathcal{C}$ is at least that
of~$\mathcal{B}$.

These ideas are clearly more general than the particular problems studied in
this paper.  However, note that we do not assert, for more general problems,
that an optimal strategy necessarily has constant regret.

We have been told that the optimality of constant-regret strategies may
be known in the context of bandit problems, although we have been unable to
find a reference for this.

\section{Optimal control for a known parameter}
\label{sec:classical-optimal}

In this section, we review the classical control problem of minimizing the
expected cost~(\ref{eq:32}) given known~$a$, via the Hamilton-Jacobi-Bellman equation~\cite{Bertsekas:2005}.

Suppose we find ourselves at position~$q$ at time~$t$.  Let $\J(q,t;a)$ be the
expected ``cost to go'':
\begin{equation*}
  \J(q,t;a) = E\bigg[\int_t^T\big(q(\tau)^2 + u(\tau)^2\big)\,d\tau\bigg],
\end{equation*}
assuming an optimal~$u$.  Then, considering a small time step $\Delta t$ and
neglecting errors small compared with~$\Delta t$, we have
\begin{equation*}
  \J(q,t;a) = \min_u \Big[(q^2 +  u^2)\Delta t + E(\J(q + \Delta q,t+\Delta t;a)\big) \Big].
\end{equation*}
Recall that $\Delta q = (a + u)\Delta t + \Delta W$, so (again neglecting errors
$o(\Delta t)$)
\begin{equation*}
  E(\Delta q) = (a+u)\,\Delta t,\qquad E\big((\Delta q)^2\big) = \Delta t.
\end{equation*}
In particular, $\Delta q$ has the order of magnitude $(\Delta t)^{1/2}$.
Consequently, Taylor expanding $\J$ to first order in~$t$ and to second order
in~$q$, we find that
\begin{subequations}
  \label{eq:opt-classical}
\begin{equation}
  \label{eq:33}
  0 = \partial_t \J + (q^2 +  u^2) + (a + u)\partial_q \J +
  \tfrac{1}{2}\partial_q^2 \J.
\end{equation}
The optimal control~$u$ minimizes the right-hand side of~(\ref{eq:33}), so
\begin{equation}
  \label{eq:34}
  u = -\frac{1}{2}\partial_q \J,
\end{equation}
and by definition,
\begin{equation}
  \label{eq:35}
  \J(q,T;a) = 0.
\end{equation}
\end{subequations}
We can guess a solution to~(\ref{eq:opt-classical}), of the form
\begin{equation}
  \label{eq:4}
  \J(q,t;a) = E_2(t) q^2 + E_1(t) qa + E_0(t) a^2 + E_\sharp(t),
\end{equation}
and obtain the following ordinary differential equations:
\begin{subequations}
\label{eq:42}
  \begin{align}
  -\dot E_2 &= 1 - E_2^2,&E_2(T) = 0 \label{eq:43}\\
  -\dot E_1 &= 2 E_2 - E_1 E_2,& E_1(T) = 0 \label{eq:44}\\
  -\dot E_0 &= E_1 - \frac{E_1^2}{4},& E_0(T) = 0\\
  -\dot E_\sharp &= E_2,&E_\sharp(T) = 0.
\end{align}
\end{subequations}
These may be solved exactly:
\begin{subequations}
\label{eq:13}
  \begin{align}
  \label{eq:6}
  E_2 &= \tanh(T - t)\\
  \label{eq:23}
  E_1 &= 2[1 - \sech(T - t)]\\
  E_0 &= (T - t) - \tanh(T - t)\\
  E_\sharp &= \log\cosh(T - t).
\end{align}
\end{subequations}
From~\eqref{eq:34} and~(\ref{eq:4}), the optimal control is then
\begin{equation}
  \label{eq:11}
  u = -E_2(t) q -\frac{E_1(t)}{2} a.
\end{equation}

For future reference, we write down the formulas analogous to~\eqref{eq:4} and~\eqref{eq:13} without the
assumption that $\lambda=1$ in equation~\eqref{eq:32}.  In place of formula~\eqref{eq:4}, we obtain
\begin{equation}
  \label{eq:26}
  \J^\lambda(q,t;a) = E_2^\lambda(t)q^2 + E_1^\lambda(t) qa + E_0^\lambda(t)a^2 + E_\sharp^\lambda(t).
\end{equation}
If we define $s=(T-t)\lambda^{-1/2}$, then in place of~\eqref{eq:13}, we obtain
\begin{subequations}
  \label{eq:20}
\begin{align}
  E_2^\lambda &= \lambda^{1/2}\tanh s\\
  E_1^\lambda &= 2\lambda(1 - \sech s)\\
  E_0^\lambda &= \lambda^{3/2}(s - \tanh s)\\
  E_\sharp^\lambda &= \lambda \log\cosh s.
\end{align}
\end{subequations}

This completes our review of the classical case, for known~$a$.

\section{A Bayesian strategy for an unknown parameter}
\label{sec:bayesian-strategy}

In this section, we solve the Bayesian control problem discussed in
Section~\ref{sec:bayesian-control}.  Now, the parameter~$a$ is unknown, but we
have a prior belief given by the probability measure~$\mu$.
Our problem exhibits an interesting feature not seen in general Bayesian control
problems.  In particular, at each time~$t$, a single real number~$\xi(t)$ captures
all the relevant information from past history up to time~$t$.  To see this we
argue as follows.

For ease of
notation, we assume for the moment that $\mu$ is given by a probability density
\begin{equation*}
  d\mu(a) = \rho(a)\,da.
\end{equation*}
We first compute the posterior probability distribution for~$a$, given history
up to time~$t$, and then use that information to find the optimal control.
To this end, we compute the joint probability for the unknown~$a$ and history up to time~$t$,
by dividing the time interval from 0 to~$t$ into small steps of duration~$\Delta
t$.  Thus we consider discrete times $\tau = 0,\Delta t, 2\Delta t,\ldots,t$.  The
joint probability density of obtaining a particular~$a$ and observing the history
$q(0),q(\Delta t),\ldots,q(t)$ up to time~$t$ is given by
\begin{equation*}
  \rho(a)\cdot\prod_{\tau=0}^{t-\Delta t} \frac{1}{\sqrt{2\pi\Delta t}}
  \exp\bigg(\frac{-\big[\Delta q(\tau) - (a + u(\tau))\Delta t\big]^2}{2\Delta t}\bigg),
\end{equation*}
because $\Delta W = \Delta q - (a + u)\Delta t$ is a normal random variable (see
equation~\eqref{eq:1}).
This joint probability density has the form
\begin{equation*}
  \rho(a) \cdot \exp\bigg(\bigg[\sum_{\tau=0}^{t-\Delta t}(\Delta q(\tau) - u(\tau)\Delta t)\bigg]a - \tfrac{1}{2}
  a^2\sum_{\tau=0}^{t - \Delta t} \Delta t\bigg)\cdot(\text{factor independent of~$a$}).
\end{equation*}
Taking the limit as $\Delta t$ tends to zero, we obtain the joint probability density
\begin{equation}
  \label{eq:2}
  \rho(a) \cdot \exp\Big(\xi a - \frac{t}{2}a^2\Big)\cdot(\text{factor
    independent of~$a$}),
\end{equation}
where
\begin{equation}
  \label{eq:3}
  \xi(t) = q(t) - q(0) - \int_0^t u(\tau)\,d\tau.
\end{equation}
Consequently, the posterior probability density for $a$, given history through
time $t$, also has the form~\eqref{eq:2}.  Passing from probability densities
$\rho$ to
general probabilities~$\mu$, we see easily that the posterior probability
measure for the unknown parameter~$a$, given history up to time~$t$, has the form
\begin{equation}
  \label{eq:39}
  d\mu_\text{posterior}(a) = d\mu(a) \cdot \exp\Big(\xi a - \frac{t}{2}a^2\Big)\cdot Z(\xi,t),
\end{equation}
where $Z(\xi,t)$ may be computed by noting that probability measures integrate to~1.

Thus, as claimed at the beginning of this section, the posterior probability measure depends on the
history $q(\tau)$ only through the single number~$\xi(t)$.  Furthermore, at a given
time~$t$, $\xi(t)$ can be computed from $q(t)$ and the history
of the control~$u$ up to time~$t$.  (Recall that we take $q(0)=0$.)  In
particular, we do not need to remember the whole history of $q(\tau)$.

We can now proceed as in Section~\ref{sec:classical-optimal}.  Suppose we find
ourselves at given values of $q$ and $\xi$ at time~$t$.  Let $\J(q,\xi,t;\mu)$ be
the expected cost to go, assuming optimal~$u$.  That is,
\begin{equation}
  \label{eq:37}
  \J(q,\xi,t;\mu) = E\bigg[\int_t^T \big(q(\tau)^2 + u(\tau)^2\big)d\tau\bigg],
\end{equation}
with $u$ picked to minimize the right hand side.  Proceeding as in
Section~\ref{sec:classical-optimal}, we may derive the following partial differential
equation for $\J$:
\begin{equation}
  \label{eq:7}
  0 = \partial_t \J + (q^2 + u^2) +  (\bar a + u)\partial_q \J + \bar a \partial_\xi \J +
  \tfrac{1}{2}\partial_q^2 \J + \partial_\xi\partial_q J +
    \tfrac{1}{2}\partial_\xi^2 J,
\end{equation}
where $\bar a(\xi,t)$ denotes the expected value of $a$ with respect to
$d\mu_\text{posterior}$, and the optimal $u$ is given by
\begin{equation}
  \label{eq:8}
  u = -\tfrac{1}{2}\partial_q \J.
\end{equation}
Again we impose the boundary condition $\J=0$ at $t=T$.

Let us specialize to the case in which our prior belief $\mu$ is a normal
distribution with mean zero and standard deviation~$\sigma$.  Then from
formula~\eqref{eq:39}, we see that our posterior belief is given by the
probability density
\begin{equation*}
  \rho(a) = \exp\Big(\xi a - \frac{a^2}{2}(t + \sigma^{-2})\Big) \cdot Z(\xi,t),
\end{equation*}
and consequently the posterior expected value of~$a$ is
\begin{equation}
  \label{eq:30}
  \bar a(\xi,t) = \frac{\xi}{t + \sigma^{-2}}.
\end{equation}
As in Section~\ref{sec:classical-optimal}, we guess a solution of the form
\begin{equation}
  \label{eq:9}
  \J(q,\xi,t) = E_2(t) q^2 + E_1(t) \bar a(\xi,t) q + \J_0(\xi,t),
\end{equation}
and obtain the following ordinary differential equations:
\begin{subequations}
\begin{align}
  \label{eq:10}
  -\dot E_2 &= 1 - E_2^2,&E_2(T) &= 0\\
  -\dot E_1 &= 2 E_2 - E_1 E_2,&E_1(T) &= 0,
\end{align}
\end{subequations}
together with a partial differential equation for $\J_0$:
\begin{equation*}
  0 = \partial_t \J_0 + E_1 \bar a^2 - \frac{1}{4} E_1^2\bar a^2 + \bar
  a\partial_\xi \J_0 + E_2 + \frac{E_1}{t + \sigma^{-2}} + \tfrac{1}{2}\partial_\xi^2 \J_0.
\end{equation*}
(It turns out that we will never need to know $\J_0$, so we do not need to solve
this partial differential equation.)
Note that the ordinary differential equations for $E_2$ and~$E_1$ are the same
as those in~\eqref{eq:42}, and thus $E_2$ and~$E_1$ are again given by
formulas~\eqref{eq:6} and~\eqref{eq:23}.  Thanks to formulas~\eqref{eq:8} and~\eqref{eq:9}, the optimal control is now given by
\begin{equation}
  \label{eq:12}
  u = -E_2(t) q - \frac{E_1(t)}{2}\bar a(\xi,t),
\end{equation}
where $\bar a$ is given by~\eqref{eq:30}.  Note that this is the same as
formula~\eqref{eq:11} for the optimal control for known $a$, but with $a$ (which
is now unknown) replaced by $\bar a(\xi,t)$.

This concludes our derivation of the Bayesian strategy.

\section{Performance of the Bayesian strategy}
\label{sec:Bayes-performance}

Now we wish to determine how well the Bayesian strategy (for unknown $a$)
performs for a particular value of~$a$.  More precisely, we now assume that $q$
evolves according to equation~\eqref{eq:1}, for a particular~$a$, where $u$ is
given by equations \eqref{eq:7},~\eqref{eq:8} for a particular prior
belief~$\mu$.
Let $\scriptJ(q,\xi,t;a,\mu)$ be the expected cost to go under the above assumptions.
As in Sections \ref{sec:classical-optimal} and~\ref{sec:bayesian-strategy}, we
can derive the following partial differential equation for $\scriptJ$:
\begin{equation}
  \label{eq:14}
  0 = \partial_t \scriptJ  + (a + u)\partial_q \scriptJ + a\partial_\xi\scriptJ
+ \tfrac{1}{2}\partial_q^2\scriptJ + \partial_\xi\partial_q\scriptJ +
\tfrac{1}{2}\partial_\xi^2\scriptJ + q^2 + u^2,
\end{equation}
with $u$ given by~\eqref{eq:7},\eqref{eq:8}.  Again, $\scriptJ=0$ when $t=T$.
Let us now specialize to the case of Gaussian prior belief, as in
Section~\ref{sec:bayesian-strategy}, so that the optimal $u$ is given by~\eqref{eq:12}.

We can guess a solution of the form
\begin{equation}
  \label{eq:15}
  \scriptJ(q,\xi,t;a,\mu) = E_2(t)q^2 + E_1(t)qa + E_0(t)a^2 + F_0(t)(\bar
  a(\xi,t)-a)^2 + F_\sharp(t).
\end{equation}
Then, thanks to~\eqref{eq:30}, we find this solution does indeed satisfy~\eqref{eq:14}, provided the following ODEs are satisfied:
\begin{subequations}
\begin{align}
  \label{eq:16}
  -\dot E_2 &= 1 - E_2^2,&E_2(T)=0\\
  -\dot E_1 &= 2 E_2 - E_1 E_2,&E_1(T)=0\\
  -\dot E_0 &= E_1 - \frac{E_1^2}{4},&E_0(T) = 0\\
  -\dot F_0 &= -\frac{2 F_0}{t+\sigma^{-2}} + \frac{E_1^2}{4},&F_0(T) = 0\label{eq:28}\\
  -\dot F_\sharp &= E_2 + \frac{F_0}{(t+\sigma^{-2})^2},&F_\sharp(T) = 0\label{eq:29}
\end{align}
\end{subequations}
Once again, $E_2$, $E_1$, and $E_0$ are as in~\eqref{eq:42}, and hence are given
by~\eqref{eq:13}.  One readily verifies that
the following satisfy equations~\eqref{eq:28} and~\eqref{eq:29}:
\begin{subequations}
  \label{eq:21}
\begin{align}
  \label{eq:17}
  F_0(t) &= (t+\sigma^{-2})^2\int_t^T \frac{E_1(\tau)^2}{4(\tau + \sigma^{-2})^2}\,d\tau,\\
  \label{eq:18}
  F_\sharp(t) &= \int_t^T \bigg[E_2(\tau) + \frac{F_0(\tau)}{(\tau + \sigma^{-2})^2}\bigg]\,d\tau.
\end{align}
\end{subequations}
This concludes our discussion of the performance of the Bayesian strategy for
given~$a$.

\section{Results for agnostic control problems}
\label{sec:results}

In this section we determine strategies that optimize each of the three variants
of agnostic control considered in Section~\ref{sec:agnostic-control}: additive
regret, multiplicative regret, and the fuel tax variant.  It turns out that in
each case, the optimal strategy is a Bayesian strategy in which the prior belief
about $a$ is a normal distribution with mean zero and standard
deviation~$\sigma$.  The optimal choice of $\sigma$ depends on which type of
agnostic control one is considering.  (Strictly speaking, for additive regret,
the optimal strategy corresponds to the limit $\sigma\to\infty$.)

\subsection{Minimizing additive regret}
\label{sec:results-AR}

In this case, as we will see in a moment, we will need to consider a
non-zero starting time~$T_0$.
The cost for the optimal control with known $a$, starting at time~$T_0$, is given by~\eqref{eq:4}:
\begin{equation*}
  \J(q,T_0;a) = E_2(T_0) q^2 + E_1(T_0) qa + E_0(T_0) a^2 + E_\sharp(T_0).
\end{equation*}
The cost for our Bayesian strategy, for a particular $a$, also starting at time
$T_0$, is given by~\eqref{eq:15}:
\begin{equation*}
  \scriptJ(q,\xi,T_0;a,\sigma) = E_2(T_0) q^2 + E_1(T_0) qa + E_0(T_0) a^2 + F_0(T_0)(\bar a(\xi,T_0)-a)^2
  + F_\sharp(T_0).
\end{equation*}
The difference between $\scriptJ$ and $\J$ is therefore
\begin{equation}
  \label{eq:38}
  F_0(T_0)(\bar a(\xi,T_0)-a)^2 + F_\sharp(T_0) - E_\sharp(T_0).
\end{equation}
Now suppose we start at time $t=0$ and position $q=0$.  Until time $T_0$ we are
required to set our control $u=0$, after which we are free to pick the
optimal~$u$.  Then $q(T_0)$ is normally distributed with
mean $aT_0$ and standard deviation $T_0^{1/2}$.  Moreover, $\xi(T_0)=q(T_0)$
(see equation~\eqref{eq:3} and recall that $q(0)=0$).  Consequently, $(\bar
a(\xi,t)-a)^2$ has expected value
\begin{equation*}
  \frac{T_0 + a^2\sigma^{-4}}{(T_0+\sigma^{-2})^2}
\end{equation*}
(see formula~\eqref{eq:30}), so
from~\eqref{eq:38} we see that the additive regret is
given by
\begin{equation}
  \label{eq:40}
  \AR(a) = F_0(T_0) \cdot\frac{T_0 + a^2\sigma^{-4}}{(T_0 + \sigma^{-2})^2} +
    F_\sharp(T_0) - E_\sharp(T_0).
\end{equation}
As $\sigma$ tends to infinity, the control strategy becomes
\begin{equation}
  \label{eq:22}
  u = -E_2(t) q - \frac{1}{2} E_1(t) \frac{\xi}{t},
\end{equation}
thanks to~\eqref{eq:30} and~\eqref{eq:12}, and the additive regret is given by
\begin{equation}
  \label{eq:41}
  F_0(T_0)T_0^{-1} + F_\sharp(T_0) - E_\sharp(T_0),
\end{equation}
which is independent of~$a$.  Thus we have found a strategy that optimizes
worst-case additive regret.  In particular, the reasoning in
Section~\ref{sec:constant-regret-opt} applies here even though our strategy is
not a Bayesian strategy for a particular prior, but is instead a limit of such
strategies.

One can check that, as $T_0$ tends to zero, the
additive regret~\eqref{eq:41} tends to infinity (in particular, $F_\sharp$
diverges logarithmically), which is why we introduced the parameter~$T_0$.
In our discussions of multiplicative regret and the fuel tax variant, we will
take $T_0=0$.

\subsection{Minimizing multiplicative regret}
\label{sec:results-MR}

\begin{figure}
  \centering
  \begin{tikzonimage}[width=0.49\textwidth]{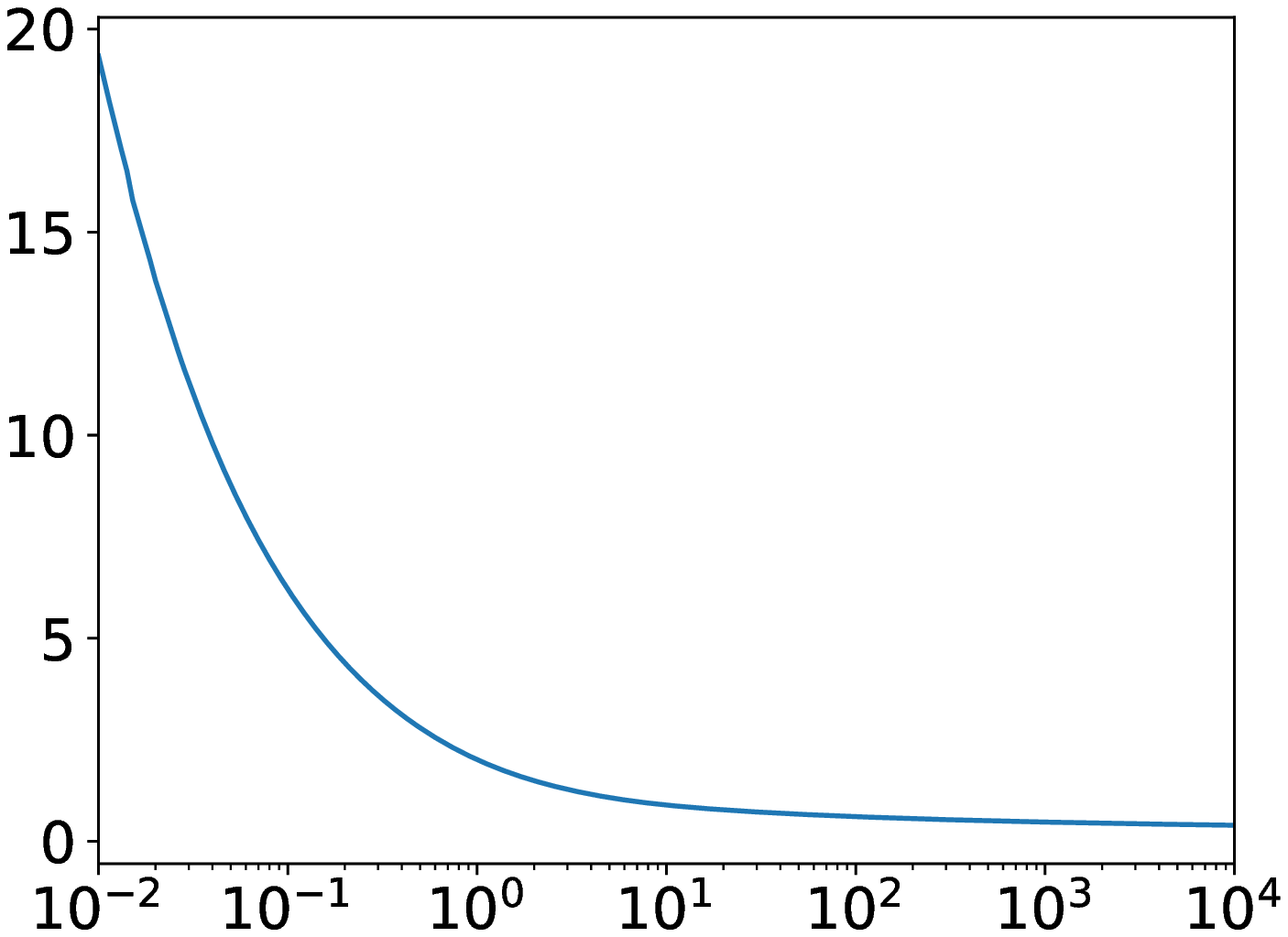}
    \node at (0.5,0.02) {$T$};
    \node at (0.02,0.5) {$\sigma$};
  \end{tikzonimage}
  \begin{tikzonimage}[width=0.49\textwidth]{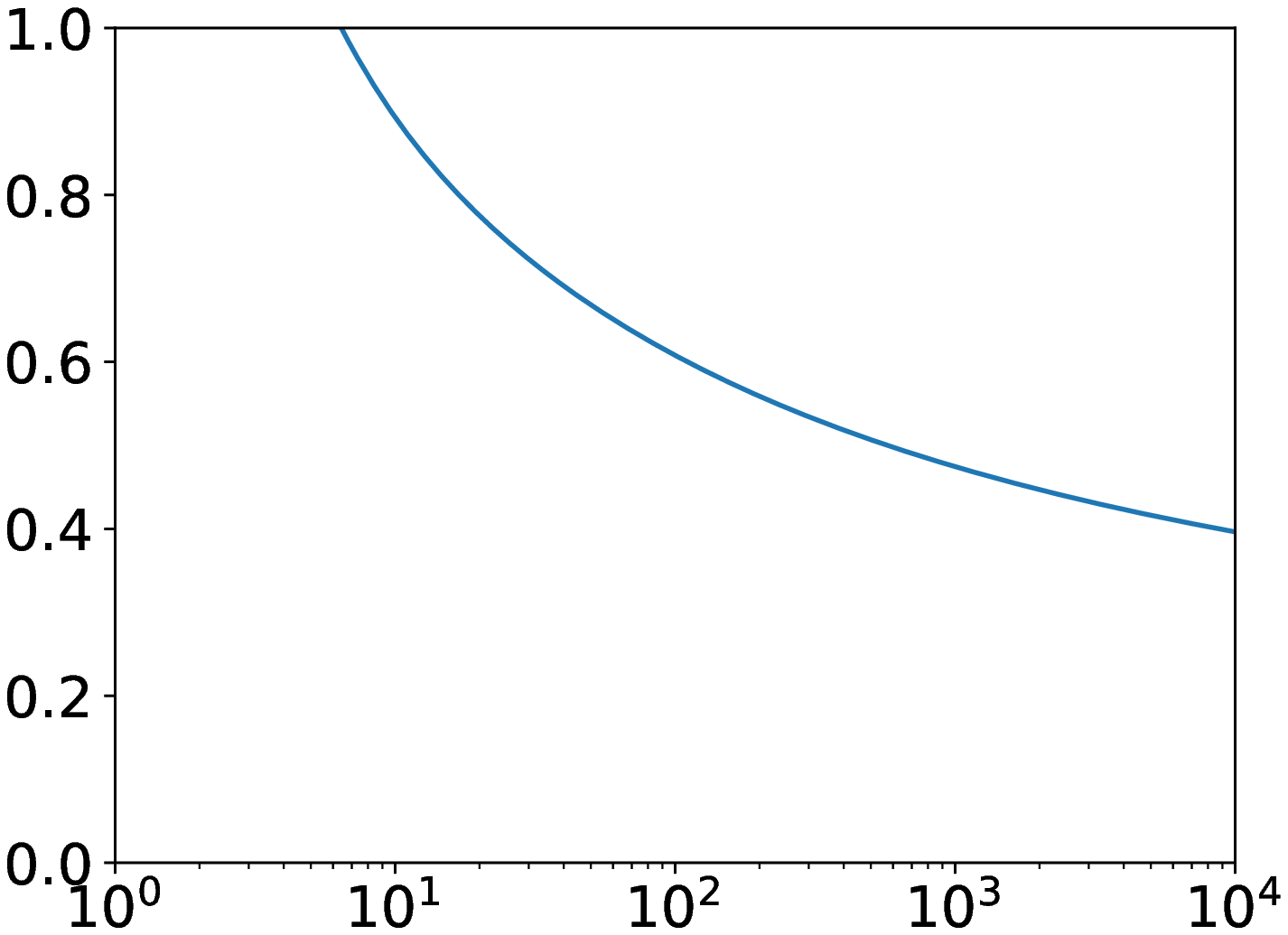}
    \node at (0.5,0.02) {$T$};
    \node at (0.02,0.5) {$\sigma$};
  \end{tikzonimage}
  \caption{Optimum standard deviation of prior belief about $a$, to minimize
    multiplicative regret.}
  \label{fig:sigma}
\end{figure}

Thanks to equations~\eqref{eq:4} and~\eqref{eq:15}, the multiplicative regret of
the optimal Bayesian strategy for normal prior belief with standard
deviation~$\sigma$ is given by
\begin{equation}
  \label{eq:19}
  \MR(a) = \frac{(E_0(0) + F_0(0))a^2 + F_\sharp(0)}{E_0(0) a^2 + E_\sharp(0)}.
\end{equation}
Here, $F_0$ and~$F_\sharp$ depend on $\sigma$; see equations~\eqref{eq:17}
and~\eqref{eq:18}.  We wish to find a value of~$\sigma$ for which the
expression~\eqref{eq:19} is independent of~$a$.  This occurs precisely when
\begin{equation}
  \label{eq:25}
  \frac{E_0(0) + F_0(0)}{E_0(0)} = \frac{F_\sharp(0)}{E_\sharp(0)}.
\end{equation}
For each~$T$, \eqref{eq:25} is a single equation for the one unknown~$\sigma$.
As a function of~$T$, the solution may be computed numerically (e.g., using
Newton's method).  To carry out the numerics, we need to evaluate the
integrals~\eqref{eq:21} using quadrature, or solve the ODEs for $F_0$ and
$F_\sharp$ backwards, from $t=T$ to $t=0$.  The optimal $\sigma$ is shown as a
function of $T$ in Figure~\ref{fig:sigma}.  One can show that the optimal
$\sigma$ tends to zero as $T\to\infty$.

\begin{figure}
  \centering
  \begin{tikzonimage}[width=0.49\textwidth]{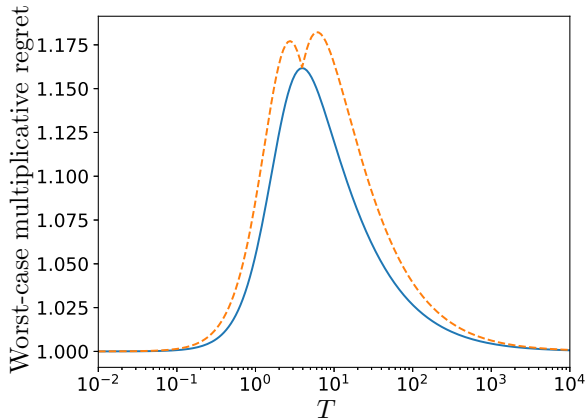}
    \node at (0.5,0.02) {$T$};
    \node[rotate=90]  at (0,0.5) {Worst-case multiplicative regret};
  \end{tikzonimage}
  \caption{Worst-case multiplicative regret for the optimum strategy (solid); and a
    Bayesian strategy with a particular $\sigma$ independent of~$T$ (dashed).}
  \label{fig:regret}
\end{figure}

Figure~\ref{fig:regret} shows the worst-case regret as a function of~$T$, for
two different strategies.  The solid curve arises from the optimum strategy: the Bayesian strategy with $\sigma$
chosen, for each $T$, to minimize worst-case multiplicative regret.  The dashed
curve arises from a Bayesian strategy, for a particular value of $\sigma$ chosen
independently of $T$.  More precisely, we pick $\sigma$ optimally for the value
of $T$ for which the solid curve reaches its maximum.
Perhaps surprisingly, both strategies perform very well.  The
optimum strategy is never more than 17\% worse than the optimum strategy for
known~$a$.  Even the ``simple'' Bayesian strategy with fixed $\sigma$ is not much
worse.  And regardless of the choice of~$\sigma$, the worst-case multiplicative regret for the Bayesian
strategy goes to $1$ as $T\to 0$ or $T\to\infty$.

\subsection{Minimizing fuel tax regret}
\label{sec:results-fueltax}

Recall from equation~\eqref{eq:32} the parameter~$\lambda$, which represents the
``cost of fuel.''  Our goal here is to compare the expected cost
$\J_\text{unknown $a$}$ of a strategy with unknown~$a$ and $\lambda=1$ with the
cost $\J_\text{known $a$}(\lambda)$ of an optimal strategy with known~$a$ and
$\lambda>1$ (see Section~\ref{sec:agnostic-control}).  As in our discussion of
multiplicative regret, we look for a Bayesian strategy with normal prior having
mean zero and standard deviation~$\sigma$.  For fixed $T$ and~$\lambda$, we pick
$\sigma$ so that the ratio $\J_\text{unknown $a$}/\J_\text{known $a$}(\lambda)$ is
independent of~$a$.  This ratio is thus a function of~$\lambda$, and we wish to
find $\lambda$ to make the ratio equal to~1.  Thanks to the discussion in
Section~\ref{sec:constant-regret-opt}, this value of $\lambda$ is then the
fuel tax regret defined in Section~\ref{sec:agnostic-control}.  Formulas \eqref{eq:26}
and~\eqref{eq:20}, together with \eqref{eq:15}, \eqref{eq:13},
and~\eqref{eq:21}, make this a routine numerical computation.

\begin{figure}
  \centering
  \begin{tikzonimage}[width=0.49\textwidth]{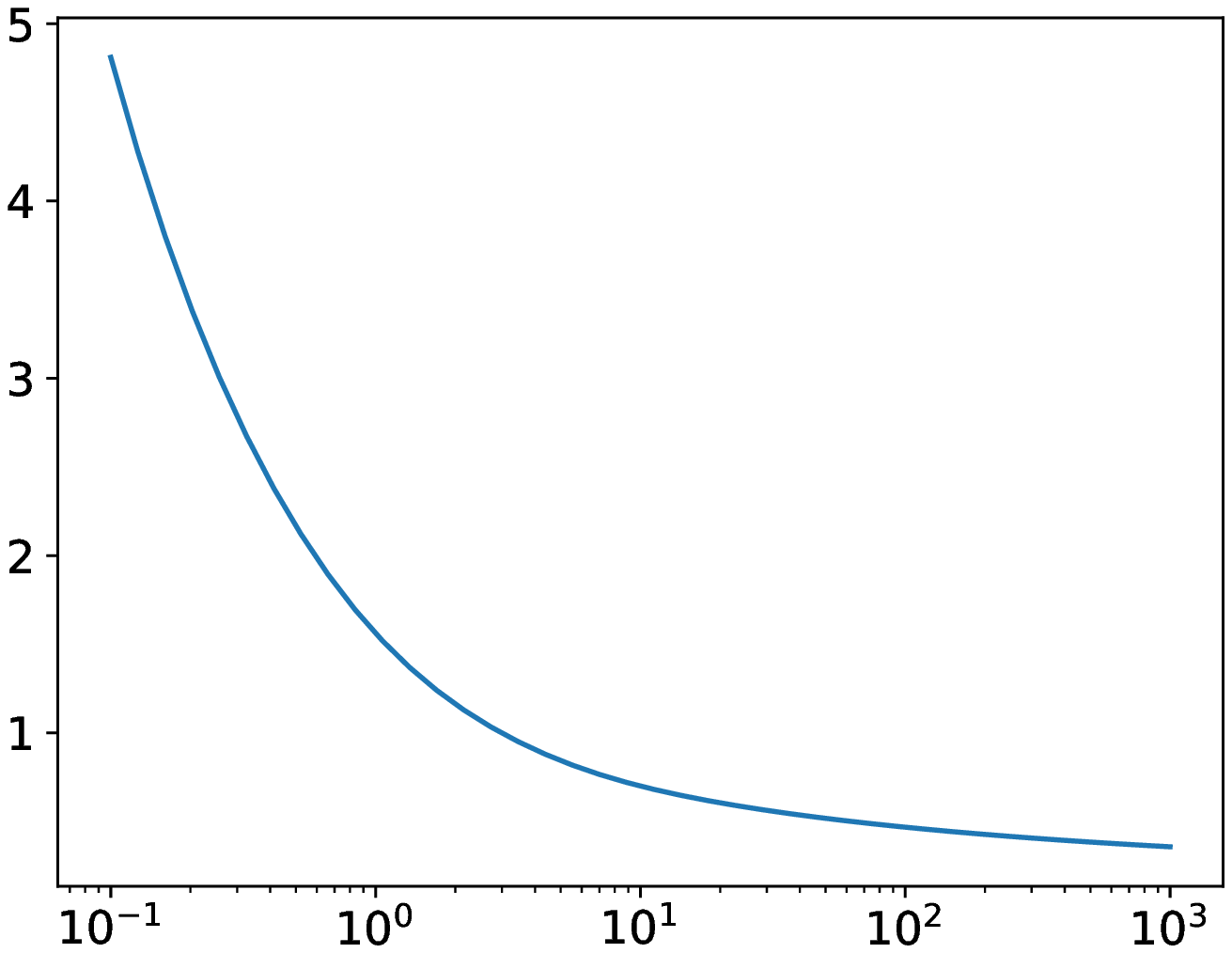}
    \node at (0.5,0.02) {$T$};
    \node at (0.02,0.5) {$\sigma$};
  \end{tikzonimage}
  \hfil
  \begin{tikzonimage}[width=0.49\textwidth]{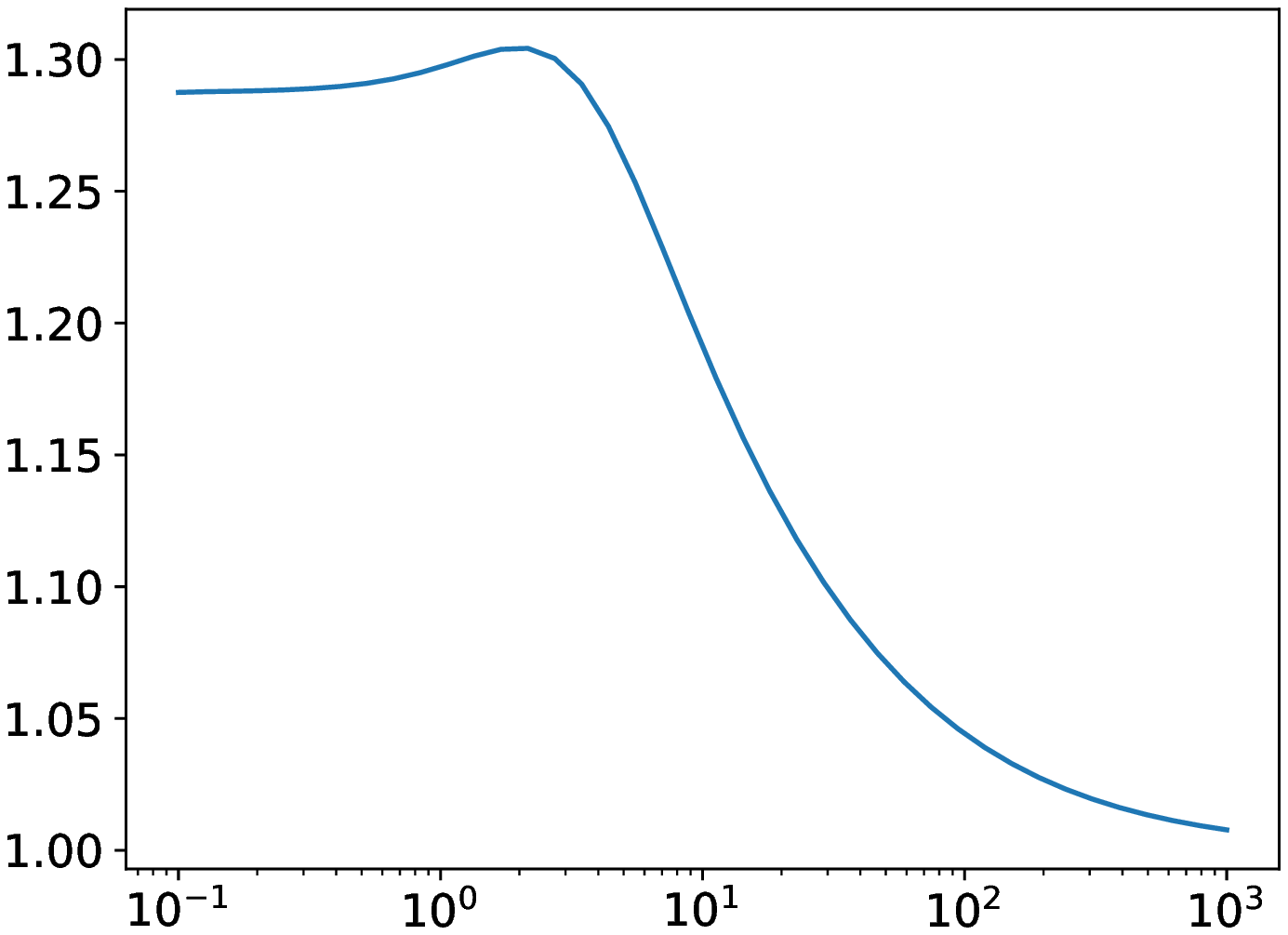}
    \node at (0.5,0.02) {$T$};
    \node at (0.02,0.5) {$\lambda$};
  \end{tikzonimage}
  \caption{Optimal Bayesian strategies for the fuel-tax variant, and resulting fuel tax regret~$\lambda$.}
  \label{fig:fuel-tax}
\end{figure}

Figure~\ref{fig:fuel-tax} shows the results of these computations.  As in the
case of multiplicative regret, the optimal value of $\sigma$ is a decreasing
function of~$T$, and we pay only a modest price for our lack of knowledge
of~$a$.  The most challenging case arises for $T\approx 2$, with a corresponding
fuel tax regret~$\lambda\approx 1.3$.

\section{Conclusions and problems for further study}
\label{sec:conclusions}
We have posed several optimal control problems for a toy model
\begin{equation*}
  \Delta q = (a + u(t)) \Delta t + \Delta W
\end{equation*}
of dynamics with noise, involving a single parameter~$a$.  These include
the classical problem in which $a$ is known, and the Bayesian problem
in which we assume a prior belief regarding~$a$.  We posed also three
``agnostic'' control problems, in which nothing is assumed about~$a$, and we
hope to minimize some notion of regret.  Here we have considered three different
notions of regret: additive regret (often called simply ``regret'');
multiplicative regret (often called ``competitive ratio''); and fuel-tax regret,
which we have not seen in the literature.

Each of these optimal control problems involves minimizing a cost of the
form
\begin{equation*}
  \J = \int_{T_0}^T (q^2 + \lambda u^2)\,dt,
\end{equation*}
for a particular terminal time~$T$.  Previous studies such as
\cite{Cohen:2019} have
produced strategies that behave well as the terminal time $T$ tends to
infinity.  We work in a different regime.  We are concerned with a fixed value
of~$T$, but assume nothing about the parameter~$a$.

For each of our agnostic control problems, minimum regret is achieved by a
Bayesian strategy (or a limit of Bayesian strategies) with the unknown $a$
assumed to be normally distributed, with mean zero and standard deviation
depending on the terminal time and the notion of regret.  The optimal strategy
for each of our agnostic control problems has regret independent of the actual
value of~$a$.

We have only begun by solving the simplest toy problem of agnostic
control.  Already, substantial challenges arise when we consider further toy
problems, with dynamics given by
\begin{equation}
  \label{eq:45}
  \Delta q = (aq + u)\Delta t + \Delta W
\end{equation}
or
\begin{equation}
  \label{eq:47}
  \Delta q = au\,\Delta t + \Delta W  
\end{equation}
in place of equation~\eqref{eq:1}.  These toy problems are the subject of
ongoing work, and we hope to return to them in a future paper.  Note that, if we
regard the unknown $a$ as part of the state, then the dynamics are linear for
our first toy problem~\eqref{eq:1}, but not for \eqref{eq:45} or~\eqref{eq:47}.
This observation may explain why toy problems \eqref{eq:45} and~\eqref{eq:47}
are more challenging than the toy problem solved here~\cite{Ramadge:PC}.

\subsection*{Acknowledgments}
This work was supported by the Air Force Office of Scientific Research, under
awards FA9550-19-1-0005 and FA9550-18-1-0069, and partially supported by the National Science Foundation,
under grant DMS-1700180.  We are particularly grateful to Fred Leve for raising the issue
of optimal control with learning on the fly, and Amir Ali Ahmadi and Elad Hazan for helpful discussions.

\bibliographystyle{siamplain}
\bibliography{regret}

\end{document}